# A RIDGE-PARAMETER APPROACH TO DECONVOLUTION

By Peter Hall and Alexander Meister

*Australian National University and Universität Stuttgart*

Kernel methods for deconvolution have attractive features, and prevail in the literature. However, they have disadvantages, which include the fact that they are usually suitable only for cases where the error distribution is infinitely supported and its characteristic function does not ever vanish. Even in these settings, optimal convergence rates are achieved by kernel estimators only when the kernel is chosen to adapt to the unknown smoothness of the target distribution. In this paper we suggest alternative ridge methods, not involving kernels in any way. We show that ridge methods (a) do not require the assumption that the error-distribution characteristic function is nonvanishing; (b) adapt themselves remarkably well to the smoothness of the target density, with the result that the degree of smoothness does not need to be directly estimated; and (c) give optimal convergence rates in a broad range of settings.

**1. Introduction.** Density estimation with observation error is almost always based on kernel methods, where the kernel depends on the error distribution. See, for example, the early contributions of Carroll and Hall [4], Liu and Taylor [22], Stefanski [27], Stefanski and Carroll [28], Zhang [30] and Fan [[12, 13, 14]]. More recent work, which also surveys earlier research, includes that of van Es, Spreij and van Zanten [29], Delaigle and Gijbels [[8, 9, 10]], Meister [[23, 24]] and Comte, Rozenholc and Taupin [[6, 7]].

However, kernel methods have disadvantages, not least being the fact that for effective implementation they require the characteristic function of the error distribution to have no zeros on the real line. In particular, the error distribution should not be compactly supported.

Motivated partly by this difficulty, in the present paper we introduce a new estimation procedure based on ridging. Since this technique does not involve a kernel, the optimal choice of which depends on the unknown smoothness of the target density, then our new method can have a relatively high degree of adaptivity. It does not require the regularity of the target function to be known in advance, and admits elementary cross-validation approaches to smoothing-parameter choice. See Fan and Koo [16] for discussion of adaptive methods in the setting of wavelet-based deconvolution.

Importantly, ridging allows us to work with error distributions that have non-positive characteristic functions. In particular, using the new method we can readily treat problems where the error distribution is compactly supported. Ridging also eliminates discontinuities in integrals, which occur (for example) when using the sinc kernel, and so avoids the need for tapers.

The ridge-based method enjoys optimal convergence rates, both in standard, or "nonoscillatory," cases where characteristic functions do not vanish, and in "oscillatory" cases where those functions have infinitely many zeros. In the latter setting, neither convergence rates nor optimality properties have been established before. The rates turn out to be particularly interesting.

For example, although as a rule optimal rates depend on the smoothnesses of both the target and the error distributions, for all sufficiently smooth target distributions, they depend only on the smoothness of the error distribution. This property is not observed in the nonoscillatory case. It makes choice of the ridge parameter remarkably straightforward; the parameter can be chosen quite easily, within a very wide range, without adversely affecting the rate of convergence.

Ridging is also adaptable to errors-in-variables problems, where it enjoys similar advantages.

**2. Methodology.** Suppose we observe data $W_1, \ldots, W_n$ generated by the model

$$W_j = X_j + \delta_j, \tag{2.1}$$

where $X_j, \delta_j$, for $1 \leq j < \infty$, are mutually independent random variables, the sequences $X_1, X_2, \ldots$ and $\delta_1, \delta_2, \ldots$ are identically distributed, and $\delta_j$ has known density $f_\delta$. We wish to estimate the density, $f_X$, say, of $X$.

Conventional estimators in this problem are given by

$$\check{f}_X(x) = \frac{1}{nh} \sum_{j=1}^{n} L\left(\frac{x - W_j}{h}\right), \tag{2.2}$$

where

$$L(u) = \frac{1}{2\pi} \int e^{-itu} \frac{K^{\mathrm{ft}}(t)}{f_\delta^{\mathrm{ft}}(t/h)} \, dt, \tag{2.3}$$

$K$ is a kernel function, $K^{\mathrm{ft}}(t) = \int e^{itx} K(x) \, dx$ is its Fourier transform, $h > 0$ is a bandwidth and we have $K^{\mathrm{ft}}(0) = 1$. Usually, $K$ is chosen so that $K^{\mathrm{ft}}$ is compactly supported, but in order for good convergence rates to be achieved, it must also take account of the unknown smoothness of the distribution of $W$. In particular, if $f_W$ has $d$ derivatives, then $K$ should be chosen so that $K^{\text{ft}}(t) = 1 + O(|t|^d)$ as $t \to 0$. See, for example, Fan [12].

Thus, effective choice of $K$ requires difficult, adaptive estimation of the smoothness of $f_X$. This problem can be alleviated by employing the sinc kernel, $K(x) = (\sin x)/(\pi x)$; then, $K^{\text{ft}}(t) = 1$ on the interval $[-1, 1]$ and vanishes elsewhere. However, in such cases the integral at (2.3) stops abruptly, before the integrand has a chance to descend to zero. That causes oscillations of Gibbs phenomenon type in the final estimator; the oscillations can usually only be removed by fitting a taper to the integrand. An approach of this type has recently been discussed by Butucea and Tsybakov [3].

A particularly significant difficulty with kernel methods arises when $f_\delta^{\text{ft}}$ vanishes at one or more points on the real line. Then there are poles in the integral at (2.3). Typically, the integrand behaves like a nonzero constant multiple of $(t-p)^{-1}$ in the neighborhood of a pole at $p$, and so the integral does not exist.

To avoid having to address this problem, it is customary in the literature to assume that $f_\delta^{\text{ft}}$ does not vanish anywhere. However, that constraint excludes all the conventional compactly supported models for the distribution of $\delta$, such as uniform and beta distributions, as well as some infinitely supported models. Devroye [11] shows that consistency is achievable whenever the set $\{t : f_\delta^{\text{ft}}(t) = 0\}$ has Lebesgue measure zero. The underlying estimator requires selection of three parameter sequences, of which one excises a neighborhood of each pole from the integral in (2.3). This technique is arguably not attractive from a practical viewpoint. Another approach, where the condition $f_\delta^{\text{ft}}(t) \neq 0$ is relaxed, is suggested by Groeneboom and Jongbloed [18]. However, their method is restricted to the case where $\delta$ has a uniform distribution.

The reason for introducing the factor $K^{\text{ft}}(t)$ to the integrand at (2.3) is to avoid difficulties when $t$ is relatively large and the denominator, $f_\delta^{\text{ft}}(t)$, is small. We suggest ridging the integrand instead. This is not completely straightforward, since the denominator in (2.3) might be negative. We propose overcoming this problem via an indirect approach, which involves first making the denominator positive and then inserting the ridge, as follows. Note that, if we multiply both the numerator and denominator at (2.3) by $f_\delta^{\text{ft}}(-t)$, we convert the denominator to $|f_\delta^{\text{ft}}(t)|^2$, a real-valued and nonnegative function. This suggests taking the integral in (2.3) over the whole real line, and incorporating a positive ridge function, $h(t)$, say, generally depending on $n$ and sometimes also on $t$.

The last step may be implemented in a variety of ways, of which one is to use

$$\tilde{f}_X(x) = \frac{1}{2\pi} \int \frac{f_\delta^{\text{ft}}(-t)|f_\delta^{\text{ft}}(t)|^r \hat{f}_W^{\text{ft}}(t)}{\{|f_\delta^{\text{ft}}(t)| \vee h(t)\}^{r+2}} e^{-itx}\, dt, \tag{2.4}$$

where $\hat{f}_W^{\text{ft}}(t) = n^{-1} \sum_j e^{itW_j}$ denotes the empirical characteristic function of $W$, $r \geq 0$ describes the "shape" of the smoothing regime, and $x \vee y$ represents the maximum of $x$ and $y$. In practice, we would confine attention to the real part of $\tilde{f}_X$, which we write as $\hat{f}_X = \Re \tilde{f}_X$. The connection between $\tilde{f}_X$, defined at (2.4), and a kernel-type estimator $\check{f}_X$, at (2.2), can be seen by taking the kernel, $L$, in (2.2) to have Fourier transform $L^{\text{ft}}$ given by

$$L^{\text{ft}}(t) = \begin{cases} h^{-(r+2)} f_\delta^{\text{ft}}(-t)|f_\delta^{\text{ft}}(t)|^r, & \text{if } |f_\delta^{\text{ft}}(t)| \leq h, \\ f_\delta(t)^{-1}, & \text{if } |f_\delta^{\text{ft}}(t)| > h. \end{cases} \tag{2.5}$$

In the nonoscillatory case, where $f_\delta^{\text{ft}}$ does not vanish on the real line, we may take $h(t)$ equal to a constant depending on $n$. In particular, optimal convergence rates, for a wide range of different smoothnesses of $f_X$, are obtained when $h(t)$ does not depend on $t$. Therefore, in such settings our approach has removed the need to choose a kernel whose shape is adapted to the unknown smoothness of $f_X$. Even when $f_\delta^{\text{ft}}$ vanishes at points on the line, it is usually straightforward to determine $h(t)$; in such cases that function is a constant multiplied by a power of $|t|$, and the power can generally be obtained from knowledge of $f_\delta$.

Note too that the approach at (2.4) removes the need for tapers, and in fact, the integrand in (2.4) is a uniformly continuous function for each $x$.

In order for (2.4) to be well defined, the integrability of $|f_\delta^{\text{ft}}|^{r+1}$ needs to be assured. This is straightforward, however. Indeed, if $f_\delta$ is square-integrable, then $r \geq 1$ is sufficient.

In "standard" errors-in-variables problems we combine the model (2.1) with a regression model, and observe data $(W_1, Y_1), \ldots, (W_n, Y_n)$, generated as

$$Y_j = g(X_j) + \epsilon_j, \qquad W_j = X_j + \delta_j, \tag{2.6}$$

where $X_j, \delta_j, \epsilon_j$, for $1 \leq j < \infty$, are mutually independent random variables, the variables $X_j$ and $\delta_j$ are as in the model (2.1) and $\epsilon_1, \epsilon_2, \ldots$ are identically distributed with $E(\epsilon_j) = 0$. We wish to estimate the smooth function $g$. There is a large literature on kernel methods in this problem; see Fan and Truong [17] for theory and Carroll, Ruppert and Stefanski [5] for a survey.

A ridge-based estimator is given by $\hat{g} = \Re \tilde{g}$, where, for any $r \geq 0$,

$$\tilde{g}(x) = \frac{\int f_\delta^{\text{ft}}(-t)|f_\delta^{\text{ft}}(t)|^r \hat{v}(t)/\{|f_\delta^{\text{ft}}(t)| \vee h(t)\}^{r+2} e^{-itx}\, dt}{\int f_\delta^{\text{ft}}(-t)|f_\delta^{\text{ft}}(t)|^r \hat{f}_W^{\text{ft}}(t)/\{|f_\delta^{\text{ft}}(t)| \vee h(t)\}^{r+2} e^{-itx}\, dt}, \tag{2.7}$$

again $\hat{f}_W^{\text{ft}}(t) = n^{-1} \sum_j e^{itW_j}$, and $\hat{v}(t) = n^{-1} \sum_j Y_j e^{itW_j}$.

In Berkson's errors-in-variables problem the observed data are $(X_1, Y_1)$, $\ldots, (X_n, Y_n)$, generated as

$$Y_j = g(X_j + \delta_j) + \epsilon_j,$$

where, as before, $X_j, \delta_j, \epsilon_j$ are mutually independent for $1 \leq j < \infty$, the $X_j$, $\delta_j$ and $\epsilon_j$ sequences are each identically distributed, $E(\epsilon_j) = 0$ and $\delta_j$ has known density $f_\delta$.

There is apparently no published account of nonparametric methods in this problem, which dates from Berkson [1] and is generally addressed using parametric or semiparametric techniques (see, e.g., Reeves et al. [26] and Buonaccorsi and Lin [2]). Using our ridge-based approach, an estimator of $g$ can be taken to be the real part of $\tilde{g}$, where

$$\tilde{g}(x) = \frac{1}{2\pi} \int \frac{f_\delta^{\text{ft}}(t)|f_\delta^{\text{ft}}(t)|^r \widehat{w}(t)}{\{|f_\delta^{\text{ft}}(t)| \vee h(t)\}^{r+2}} e^{-itx}\, dt,$$

$\widehat{w}(t) = \sum_j D_j Y_j e^{itX_j}$ and $D_j > 0$ denotes the distance from $X_j$ to the nearest other data value. The function $\widehat{w}(t)$ estimates $c^{\text{ft}}(t)/f_\delta^{\text{ft}}(t)$, where $c^{\text{ft}}$ is the Fourier transform of the function $c(x) = \int f_\delta(x-u) E(Y|X=u)\, du$.

**3. Smoothing-parameter choice.** When $f_\delta^{\text{ft}}$ does not vanish on the real line and also in the case of supersmooth $f_\delta$, where $f_\delta^{\text{ft}}$ decreases exponentially fast to zero in the tails, it is usually adequate to take $h = h(t)$ in (2.4) to be a constant depending on $n$. In this case, $h$ is the single smoothing parameter on which the methodology depends. In contexts where $f_\delta$ is ordinary-smooth (i.e., $f_\delta^{\text{ft}}$ decreases only polynomially fast) and is compactly supported, we usually need to take $h$ to be a polynomial in $t$. Taking these cases together, we might consider

(3.1) $$h(t) = h_n(t) = n^{-\zeta}|t|^\rho,$$

where $\zeta > 0$ and $\rho \geq 0$.

Section 4 will discuss choices of $\zeta$ and $\rho$ that lead to optimal rates of convergence. The case where $f_\delta$ is compactly supported is particularly interesting; we treat it here through an example, as follows. If $f_\delta$ is the $\mu$-fold convolution of a symmetric uniform density, where $\mu \geq 2$, and if $f_X$ has an integrable second derivative, then optimal convergence rates are achieved with $\rho = 2$ and $\zeta = \frac{1}{2}$ in (3.1).

In this setting we might take $\rho = 2$ and choose the constant $\xi$ in the ridge-parameter formula $h(t) = \xi t^2$ empirically; in the context where $f_\delta^{\text{ft}}$ does not vanish, we can take $\rho = 0$ and choose $\xi$ in $h(t) = \xi$ empirically. Therefore, we should address the case where $h(t) \equiv \xi|t|^\rho$, with $\rho \geq 0$ known, and discuss selection of $\xi$. Our approach to solving this problem will be via cross-validation; see Hesse [20] for an account of this method in the case of kernel-based deconvolution.

Let $\hat{f}_X(\xi, x)$ denote the density estimator at (2.4) for this choice of the ridge. The aim is to minimize

$$\int |\hat{f}_X(\xi; x) - f_X(x)|^2 \, dx,$$

or equivalently, using the Plancherel identity, to minimize the function $J(\xi) - 2\Re I(\xi)$, where $J(\xi) = \int |\hat{f}_X(\xi; x)|^2 \, dx$ and

$$I(\xi) = \int \hat{f}_X(\xi; x) f_X(x) \, dx = \frac{1}{2\pi} \int \frac{|f_\delta^{\mathrm{ft}}(t)|^r}{\{|f_\delta^{\mathrm{ft}}(t)| \vee \xi |t|^\rho\}^{r+2}} \hat{f}_W^{\mathrm{ft}}(t) f_W^{\mathrm{ft}}(-t) \, dt.$$

While $J(\xi)$ is known, we have to produce an accessible version of $I(\xi)$ to eliminate the unknown $f_W^{\mathrm{ft}}$. That is given by

$$\hat{I}(\xi) = \frac{1}{2\pi n(n-1)} \int \frac{|f_\delta^{\mathrm{ft}}(t)|^r}{\{|f_\delta^{\mathrm{ft}}(t)| \vee \xi |t|^\rho\}^{r+2}} \sum\sum_{j \neq k} \exp\{it(W_j - W_k)\} \, dt,$$

and so we use as our criterion

(3.2) $$\mathrm{CV}(\xi) = J(\xi) - 2\Re \hat{I}(\xi).$$

(Note that $r > 0$ has to be chosen sufficiently large to ensure the integrability of $|f_\delta^{\mathrm{ft}}|^r$, and that, provided $f_\delta$ is square-integrable, this requires only $r \geq 2$.) Finally, we select the smoothing parameter

(3.3) $$\hat{\xi} = \arg\min_{\xi > 0} \mathrm{CV}(\xi).$$

Sections 4 and 5 will briefly discuss theoretical and numerical properties, respectively, of this technique.

This approach has a straightforward analogue for determining smoothing parameters in the "standard" errors-in-variables problem (2.6). In that case we treat separately the numerator and denominator in (2.7), and so the estimator $\tilde{g}$ is computed using two different smoothing parameters. The case of Berkson's errors-in-variables problem is more difficult to address, however.

## 4. Theoretical properties.

4.1. *The nonoscillatory case.* Here we consider the case where the characteristic function of $\delta$ does not vanish on the real line. Given $\beta > 1/2, C > 0$, define $\mathcal{F}_{\beta C}^1$ to be the Sobolev class of all densities $f_X$ for which

(4.1) $$\int |f_X^{\mathrm{ft}}(t)|^2 (1 + t^2)^\beta \, dt \leq C.$$

Concerning the error density, $f_\delta$, we consider ordinary-smooth densities

(4.2) $\quad C_1(1+t^2)^{-\nu} \leq |f_\delta^{\text{ft}}(t)|^2 \leq C_2(1+t^2)^{-\nu}, \qquad \text{for } -\infty < t < \infty,$

with $\nu > 0$ and $0 < C_1 < C_2 < \infty$, as well as supersmooth densities

(4.3) $\quad \exp(-c_1|t|^\gamma) \leq |f_\delta^{\text{ft}}(t)| \leq \exp(-c_2|t|^\gamma), \qquad \text{for } -\infty < t < \infty,$

with $c_1 \geq c_2 > 0$ and $\gamma > 0$.

We take the ridge function $h(t)$ to be a scalar, that is, $\rho = 0$ in the notation of (3.1). Let $E_{X\delta}$ denote expectation under the assumption that $X_1, \ldots, X_n$ and $\delta_1, \ldots, \delta_n$ have densities $f_X$ and $f_\delta$, respectively. Write $|\cdot|$ for the $L_2$ norm on the space of square-integrable, real-valued functions.

THEOREM 4.1. *In* (3.1) *take* $\rho = 0$ *and* $\zeta$ *as given below.* (a) *If* $f_\delta$ *satisfies* (4.2), *and if* $r > 0 \vee (\nu^{-1} - 1)$ *and* $\zeta = \nu/(2\beta + 2\nu + 1)$, *then the estimator* $\hat{f}_X = \Re \tilde{f}_X$, *with* $\tilde{f}_X$ *defined at* (2.4), *satisfies*

$$\sup_{f_X \in \mathcal{F}_{\beta C}^1} E_{X\delta} \|\hat{f}_X - f_X\|^2 = O(n^{-2\beta/(2\beta + 2\nu + 1)}).$$

(b) *If* $f_\delta$ *satisfies* (4.3), *and if* $r = 0$ *and* $0 < \zeta < \frac{1}{4}$, *then*

$$\sup_{f_X \in \mathcal{F}_{\beta C}^1} E_{X\delta} \|\hat{f}_X - f_X\|^2 = O\{(\log n)^{-2\beta/\gamma}\}.$$

Optimality of these convergence rates follows from results of Fan [[12, 15]] for Hölder classes. Under additional regularity assumptions on $f_\delta$, they can be extended to Sobolev classes; see Neumann [25] and Hesse and Meister [21]. A proof of Theorem 4.1 is included in a longer version of this paper (Hall and Meister [19])

4.2. *The oscillatory case.* Given $\beta > 1/2, C > 0$, define $\mathcal{F}_{\beta C}^2$ and $\mathcal{F}_{\beta C}^3$ to be, respectively, the classes of densities $f_X$ for which

$$\int \{|f_X^{\text{ft}}(t)|^2 + |f_X^{\text{ft}}(t)||(f_X^{\text{ft}})'(t)|\}(1+t^2)^\beta \, dt \leq C$$

and $|f_X^{\text{ft}}(t)| \leq C|t|^{-\beta-(1/2)}$. [In the case of $\mathcal{F}_{\beta C}^2$ we assume that $(f_X^{\text{ft}})'$ is well defined.] These conditions amount to upper bounds on the smoothness of $f_X$, alternative to that given by (4.1). To appreciate the connections among $\mathcal{F}_{\beta C}^1$, $\mathcal{F}_{\beta C}^2$ and $\mathcal{F}_{\beta C}^3$, note that $\mathcal{F}_{\beta C}^2 \subseteq \mathcal{F}_{\beta C}^1$, and that, provided

(4.4) $\qquad\qquad\qquad |(f_X^{\text{ft}})'|/|f_X^{\text{ft}}| \leq C_1$

(which condition is typically true for Laplace-type distributions, for which $|f_X^{\text{ft}}(t)|$ decreases in a polynomial way as $|t|$ increases), $f_X \in \mathcal{F}_{\beta C}^1$ entails $f_X \in \mathcal{F}^2_{\beta C_2}$, where $C_2 = C(1 + C_1)$. Also, if $f_X \in \mathcal{F}^3_{\beta C_3}$ and $\epsilon \in (0, \beta - \frac{1}{2})$, then $f_X \in \mathcal{F}^2_{\beta-\epsilon, C_3}$, where

$$C_3 = \int_0^1 (1+t^2)^\beta \, dt + C^2 \int_1^\infty |t|^{-2\beta-1}(1+t^2)^\beta \, dt;$$

and if $f_X \in \mathcal{F}^3_{\beta C}$ and (4.4) holds, then $f_X \in \mathcal{F}^1_{\beta-\epsilon, C_4}$, where

$$C_4 = (C_1 + 1) \int_0^1 (1+t^2)^\beta \, dt + C(C_1 + C) \int_1^\infty |t|^{-2\beta-1}(1+t^2)^\beta \, dt.$$

To give a more intuitive description of the densities in $\mathcal{F}^3_{\beta C}$, we mention that, for integer $\beta + \frac{1}{2}$, the relation $f_X \in \mathcal{F}^3_{\beta C}$ follows if the derivatives $f_X^{(l)}(x)$ tend to zero as $|t| \to \infty$, for all $l \leq \beta + \frac{1}{2}$, and $\int |f^{(\beta+1/2)}(x)| \, dx \leq C$. Hence, $\mathcal{F}^3_{\beta C}$ might be interpreted as an $L_1(\mathbb{R})$-analogue of the Sobolev class $\mathbb{F}^1_{\beta C}$.

Given $\mu \geq 1$, $\nu > 0$, $0 < C_1 < C_2 < \infty$, $\lambda > 0$ and $T > 0$, denote by $\mathcal{G}_{\nu\mu\lambda}$ the class of probability densities $f_\delta$ for which

(4.5) $\quad C_1 |\sin(\lambda t)|^\mu |t|^{-\nu} \leq |f_\delta^{\text{ft}}(t)| \leq C_2 |\sin(\lambda t)|^\mu |t|^{-\nu} \qquad$ for all $|t| > T$

and $f_\delta^{\text{ft}}(t)$ does not vanish for $|t| \leq T$.

The parameter $\mu$ describes the "order" of the isolated zeros of $f_\delta^{\text{ft}}$. Note that all self-convolved uniform densities are in $\mathcal{G}_{\nu\mu\lambda}$ for appropriate parameter choices, as too are their convolutions with any ordinary-smooth density. Accordingly, we introduce the class $\mathcal{G}'_{d\gamma\mu\lambda}$ of all densities satisfying (4.5) when $|t|^{-\nu}$ is replaced by $\exp(-d|t|^\gamma)$, with $d, \gamma > 0$. For instance, convolutions of uniform densities with normal densities are included in $\mathcal{G}'_{d\gamma\mu\lambda}$ for suitably chosen parameters.

Preparing for Theorem 4.2(a), put

(4.6) $\quad \rho \begin{cases} \in \left(\dfrac{\mu+\nu}{2\mu-1}, 2\mu\beta - \nu\right), & \text{if } 2\beta + 2\nu + 1 < 4\mu\beta, \\ = \dfrac{\mu+\nu}{2\mu-1}, & \text{if } 2\beta + 2\nu + 1 = 4\mu\beta, \\ \in \left(2\mu\beta - \nu, \dfrac{\mu+\nu}{2\mu-1}\right), & \text{otherwise,} \end{cases}$

(4.7) $\quad \zeta = \begin{cases} \dfrac{1}{2}, & \text{if } 2\beta + 2\nu + 1 \leq 4\mu\beta, \\ \dfrac{\nu+\rho}{2\beta+2\nu+1}, & \text{otherwise.} \end{cases}$

THEOREM 4.2. (a) *Define $h(t)$ as at (3.1), with $\rho$ and $\zeta$ as in (4.6) and (4.7), respectively. If $f_\delta \in \mathcal{G}_{\nu\mu\lambda}$ where $\nu > 0$ and $\mu \geq 1$, and if $r > 0 \vee (\nu^{-1} -$*

1), *then, for* $\hat{f}_X = \Re \tilde{f}_X$ *and* $j = 2, 3$, *we have*

$$\sup_{f_X \in \mathcal{F}^j_{\beta C}} E_{X\delta} \|\hat{f}_X - f_X\|^2$$

(4.8)
$$= \begin{cases} O(n^{-1/(2\mu)}), & \text{if } 2\beta + 2\nu + 1 < 4\beta\mu, \\ O(n^{-1/(2\mu)} \log n), & \text{if } 2\beta + 2\nu + 1 = 4\beta\mu, \\ O(n^{-2\beta/(2\beta+2\nu+1)}), & \text{otherwise.} \end{cases}$$

(b) *Define* $h(t)$ *as at* (3.1), *with* $\rho = 0$ *and* $\zeta \in (0, \frac{1}{4})$. *If* $f_\delta \in \mathcal{G}'_{d\gamma\mu\lambda}$, *and if* $r = 0$, *then for* $j = 1, 2, 3$,

$$\sup_{f_X \in \mathcal{F}^j_{\beta C}} E_{X\delta} \|\hat{f}_X - f_X\|^2 = O\{(\log n)^{-2\beta/\gamma}\}.$$

To interpret part (a) of the theorem, note that for ordinary-smooth $f_\delta$, the non-classical rate $n^{-1/(2\mu)}$ arises if $\mu$ is large enough in relation to $\nu$ and $\beta$. Then the difficulty created by the isolated zeros of $f_\delta^{\text{ft}}$ dominates the difficulty caused by the decay of $f_\delta^{\text{ft}}$ in the tails. Part (b) implies that, unlike the case of ordinary-smooth densities, the existence of isolated zeros of $f_\delta^{\text{ft}}$ does not cause any deterioration of convergence rates for supersmooth error densities, even for the comprehensive smoothness class $\mathcal{F}^1_{\beta C}$.

Next we show that the convergence rates in Theorem 4.2(a) are optimal, or nearly optimal in the case $2\beta + 2\nu + 1 = 4\beta\mu$.

THEOREM 4.3. *Let* $\hat{f}$ *be an arbitrary estimator of* $f$ *based on the sample* $W_1, \ldots, W_n$. *Assume the existence of at least one* $T > 0$ *such that*

(4.9)
$$\limsup_{t \to T} |f_\delta^{\text{ft}}(t)| / |t - T|^\mu < \infty \quad \text{and}$$
$$\limsup_{t \to T} |(f_\delta^{\text{ft}})'(t)| / |t - T|^{\mu-1} < \infty.$$

*Then for* $j = 2, 3$,

(4.10)
$$\limsup_{n \to \infty} \sup_{f_X \in \mathcal{F}^j_{\beta C}} n^{1/(2\mu)} E_{X\delta} \|\hat{f} - f_X\|^2 > 0.$$

*If, in addition,* $f_\delta \in \mathcal{G}_{\nu\mu\lambda}$ *and* $|(f_\delta^{\text{ft}})'(t)| \leq c|t|^{-\nu}$ *for some* $c > 0$ *and all* $t$, *and* $C$ *is sufficiently large, then*

(4.11)
$$\limsup_{n \to \infty} \sup_{f_X \in \mathcal{F}^3_{\beta C}} n^{2\beta/(2\beta+2\nu+1)} E_{X\delta} \|\hat{f} - f_X\|^2 > 0.$$

Note that, for integer $\mu$, (4.9) is satisfied if $f_\delta^{\text{ft}}$ is $\mu$-times continuously differentiable in a neighborhood of $T$, and $(f_\delta^{\text{ft}})^{(0)}(T) = \cdots = (f_\delta^{\text{ft}})^{(\mu-1)}(T) =$

0. Due to the greater stringency of the smoothness classes when $j = 2, 3$, the lower bound (4.11) does not follow from earlier results for $j = 1$.

Note that (4.9), in contradistinction to (4.5), requires $f_\delta^{\text{ft}}$ to decrease to zero at only a single point. Hence, if $f_\delta^{\text{ft}}$ decreases to zero at different polynomial rates at different points, and if the fastest of these rates is $\kappa$, say, then the convergence rate of an arbitrary estimator $\hat{f}$ to $f_X$ can be no faster than $n^{-1/(2\kappa)}$. Also note that (4.10) remains valid if we replace $f_X \in \mathcal{F}_{\beta C}^j$ by stronger smoothness assumptions, for instance, classes of supersmooth $f_X$ where $f_X^{\text{ft}}$ shows exponential decay or densities whose Fourier transforms are compactly supported as long as the endpoint of their support is larger than $T$.

Assumption (4.9) is weaker in other respects than (4.5), in particular, with regard to tail behavior of $f_\delta^{\text{ft}}$. Although the second part of (4.9) involves an assumption on the derivative of $f_\delta^{\text{ft}}$, that condition is a natural reflection of the first part of (4.9).

We mention that slower minimax results are derived for the estimator (2.4) under classes of ordinary smooth densities which are weaker than $\mathcal{F}_{\beta C}^2$ or $\mathcal{F}_{\beta C}^3$; see the long version of this paper for details (Hall and Meister [19]).

4.3. *Adaptivity.* Here we state an optimality result for data-driven selection of the smoothing parameter, $\hat{\xi}$, defined at (3.3). For simplicity, we confine attention to the nonoscillatory case and a special oscillatory case.

THEOREM 4.4. *Assume either the conditions of Theorem* 4.1(a) *with* $f_X \in \mathcal{F}_{\beta C}^1$, *or the conditions of Theorem* 4.2(a) *with* $f_X \in \mathcal{F}_{\beta C}^2 \cup \mathcal{F}_{\beta C}^3$, *and* $r \geq 2 \vee (2/\nu)$, $4\mu\beta < \beta + 2\nu + 1$ *and* $\rho < (2\mu + \nu)/(4\mu - 1)$. *Then, with probability* 1,

$$\text{(4.12)} \qquad \frac{\|\hat{f}_{\hat{\xi}} - f_X\|^2}{\inf_{\xi>0} E\|\hat{f}_\xi - f_X\|^2} \to 1.$$

Note that the parameters $r$ and $\rho$ do not depend on $\beta$, and give the "scale" of smoothness classes to which the choice of $\xi$ adapts.

A proof of Theorem 4.4 is included in the long version of this paper (Hall and Meister [19]).

In related work, although in the setting of kernel methods, Delaigle and Gijbels [9] compare plug-in and bootstrap methods for choosing the bandwidth, and Delaigle and Gijbels [10] suggest a bootstrap technique. The approaches discussed in both articles produce, like our cross-validation algorithm, asymptotic optimality. A major difference, however, is that in our setting, for our nonkernel method, the level of smoothness is not supposed known in advance. In this context, Theorem 4.4 shows that cross-validation can choose the degree of smoothness adaptively. By way of contrast, the selection of the kernel, in the case of a kernel method, is made in the light of the assumed level of smoothness of the target density, and theoretical arguments are predicated on that level being correct.

**5. Numerical properties.** Here we present simulation results addressing the performance of our estimators. We give graphs in two cases, (a) the regular-smooth, nonoscillatory case, and (b) the regular-smooth, oscillatory case.

In each setting, and for each parameter setting, we drew 100 random samples and ranked them in order of the size of integrated squared error (ISE). In each figure the five unbroken curves are density estimates corresponding to the largest, 25th largest, 50th largest, 75th largest and smallest value of ISE; the dashed curve depicts the true density, $f_X$. Sample size, $n$, is given below each figure. As discussed in Sections 3 and 4, in cases (a) and (b) we fixed $\rho$ at 0 and 2, respectively, and chose the only remaining smoothing parameter, $\xi$ in the formula $h(t) = \xi |t|^\rho$, by cross-validation.

Numerical results in case (a), but using methods quite different from our own, are widely available in the literature. For recent examples, see Delaigle and Gijbels [[9, 10]] in the setting of kernel techniques, and Comte, Rozenholc and Taupin [[6, 7]] for penalization methods. Case (a) is also addressed in Figures 1–6, where we take $f_\delta$ to be the Laplace density. As the target density $f_X$ we used the bimodal $\frac{1}{2}\{N(2,1) + N(-2,1)\}$ density for Figures 1 and 2, the two-fold convolution of the Laplace density for Figures 3 and 4, and the shifted $\chi^2$-density $f_X(x) = (1/16)(x+4)^2 \exp\{-\frac{1}{2}(x+4)\}$, for $x > -4$, for Figures 5 and 6.

In Figures 7–12 we illustrate case (b), for the same respective densities $f_X$ as in case (a). In case (b), $f_\delta$ was the uniform density on $[-1, 1]$.

It can be deduced from the figures for cases (a) and (b) that: (i) our estimator has a little more difficulty with the bimodal density, compared with the unimodal ones; (ii) the estimator finds the nonoscillatory case (a) more challenging than the oscillatory case (b); and (iii) performance gradually improves with increasing sample size. Property (i) is to be expected, since case (a) is characterized by greater structure, which the estimator is pressed to discover; property (ii) is the result of $f_\delta^{\text{ft}}$ decreasing more rapidly in the tails in case (a) than in case (b); and property (iii) reflects a steady decline in values of averaged integrated squared error as $n$ grows. Simulations in the case of supersmooth error, more precisely, $f_\delta = N(0,1)$ and $f_\delta$ equal to the convolution of a normal and a uniform density, show that $n = 1000$ gives results broadly similar to those for $n = 400$ in case (b).

Figures 13 and 14 show results for the kernel estimator in the setting of Figures 1 and 2. Sample sizes are $n = 400$ and 700 for the respective figures, the density in each case is the normal mixture defined three paragraphs above, and the smoothing parameter (this time, the bandwidth) was chosen by cross-validation. The kernel for these results has Fourier transform $(1 - t^2)^3$ for $|t| \leq 1$; this choice is popular in kernel deconvolution. For both sample sizes the kernel estimator is more erratic than its ridge competitor, reflecting the fact that it has consistently higher values of average mean integrated squared error. This is observed in the case of the unimodal density too.

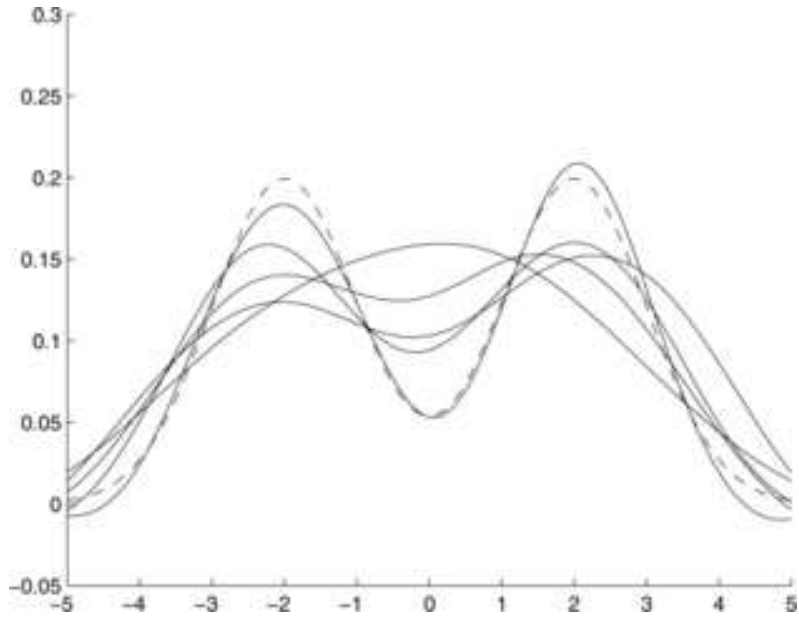

Fig. 1. $n = 400$, $r = 2$, $\rho = 0$, $\xi$ by CV average integrated squared error $= 0.013$.

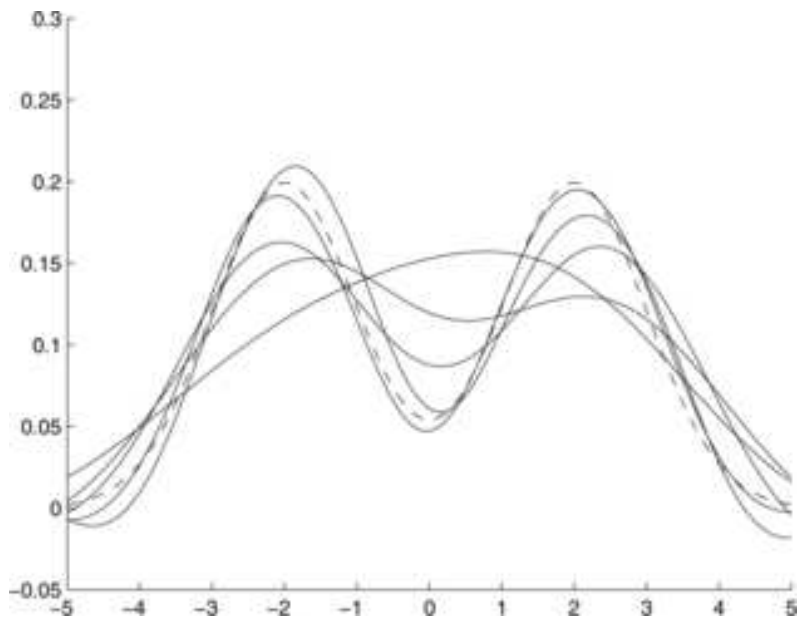

Fig. 2. $n = 700$, $r = 2$, $\rho = 0$, $\xi$ by CV average integrated squared error $= 0.0009$.

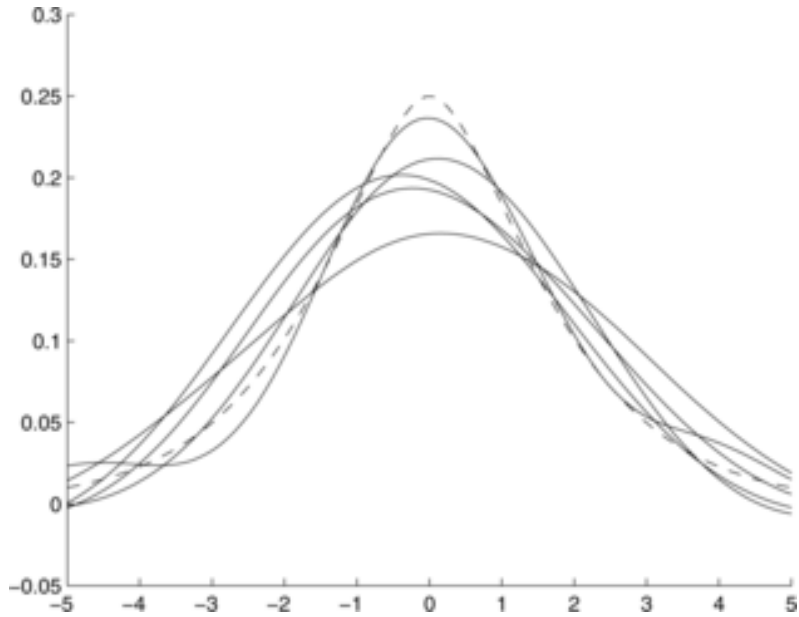

FIG. 3. $n = 400$, $r = 2$, $\rho = 0$, $\xi$ by CV average integrated squared error $= 0.006$.

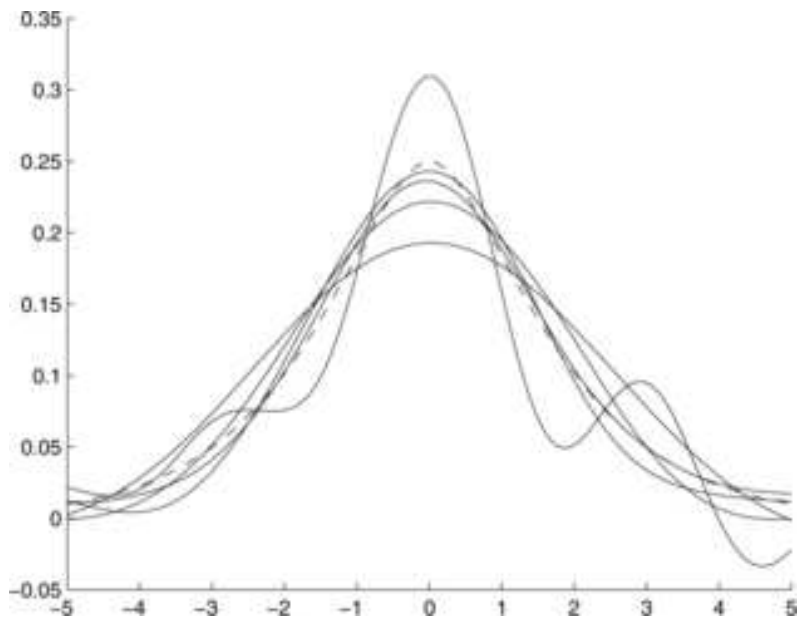

FIG. 4. $n = 700$, $r = 2$, $\rho = 0$, $\xi$ by CV average integrated squared error $= 0.004$.

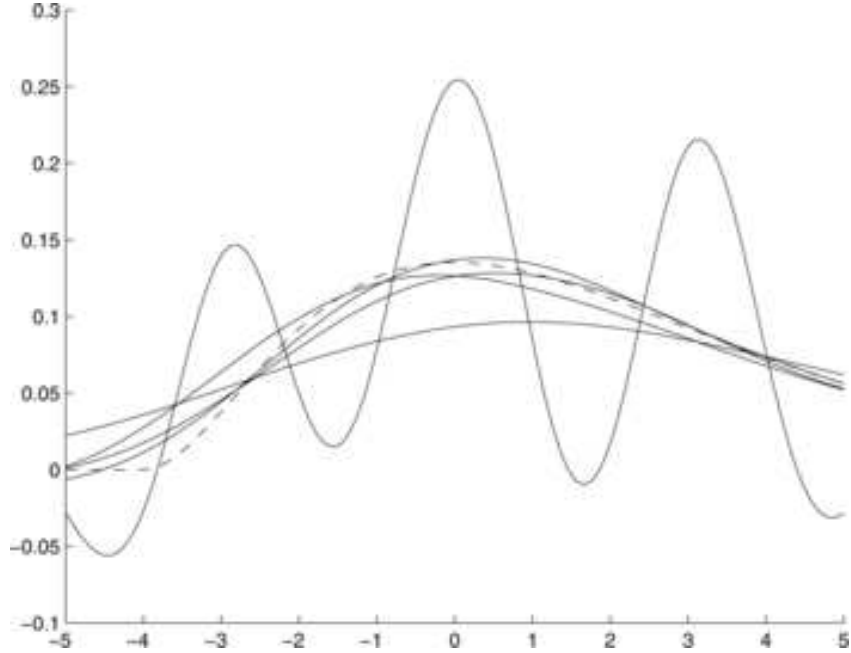

Fig. 5. $n = 400$, $r = 2$, $\rho = 0$, $\xi$ by CV average integrated squared error $= 0.005$.

**6. Outline proofs.** In this section we write const. for a generic positive constant.

6.1. *Preparatory lemma.* The following result can be proved by Parseval's identity and Fubini's theorem. Define $G_n = \{t \in \mathbb{R} : |f_\delta^{\text{ft}}(t)|^2 < h(t)^2\}$, and write $G_n^c$ for the complement of $G_n$.

LEMMA 6.1. *The mean integrated squared error of our estimator is bounded above by $V_n + B_n$ and by $V_{1,n} + V_{2,n} + B_n$, where*

$$V_n = \frac{1}{2\pi} n^{-1} \int |f_\delta^{\text{ft}}(t)|^{2+2r} h(t)^{-(4+2r)} \, dt,$$

$$V_{1,n} = \frac{1}{2\pi} n^{-1} \int_{G_n} |f_\delta^{\text{ft}}(t)|^{2+2r} h(t)^{-(4+2r)} \, dt,$$

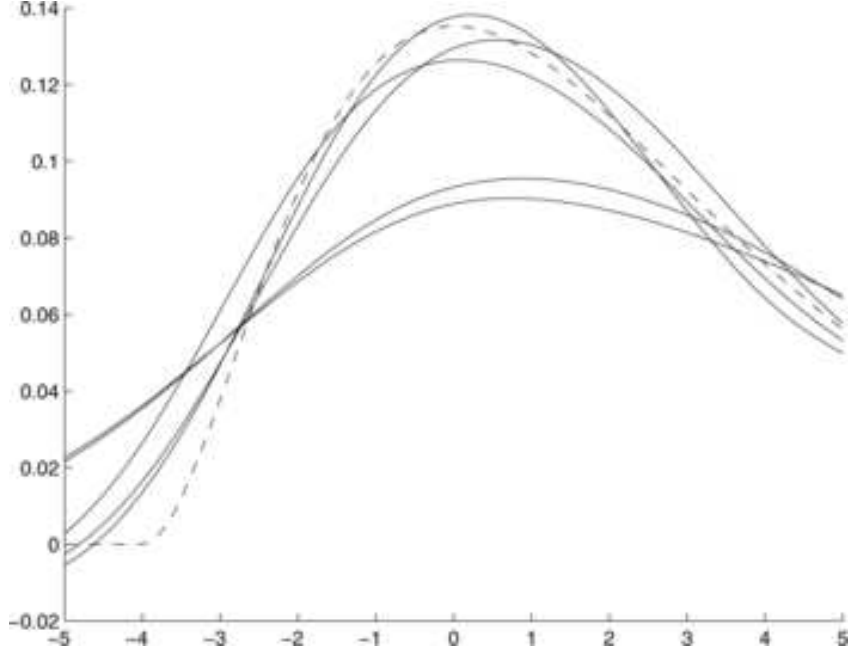

Fig. 6. $n = 700$, $r = 2$, $\rho = 0$, $\xi$ by CV average integrated squared error $= 0.004$.

$$V_{2,n} = \frac{1}{2\pi} n^{-1} \int_{G_n^c} |f_\delta^{\text{ft}}(t)|^{-2} \, dt, \qquad B_n = \frac{1}{2\pi} \sup_{f_X} \int_{G_n} |f_X^{\text{ft}}(t)|^2 \, dt.$$

The bounds $V_n + B_n$ and $V_{1,n} + V_{2,n} + B_n$ will be used to derive rates for supersmooth and ordinary-smooth error densities, respectively. In particular, Lemma 6.1 leads directly to Theorem 4.1.

6.2. *Proof of Theorem* 4.2(a). Define $G_{n,j} = G_n \cap I_j$ and $I_j = [t_{j,-}, t_{j,+}]$, where $t_{j,\pm} = (j \pm \frac{1}{2})\pi/\lambda$. Since $G_n$ is symmetric about zero, and $G_{n,0}$ is empty for $n$ sufficiently large, we may restrict attention to $j \geq 1$. Now, for $t \in G_{n,j}$,

$$|\sin(\lambda t)| \leq C_1^{-1/\mu} n^{-\zeta/\mu} |t|^{(\rho+\nu)/\mu} = \varphi_{1,n}(t),$$

say. It may be shown using a geometric argument that $G_{n,j} \subseteq [t'_{n,j,-}, t'_{n,j,+}]$, where, taking the plus and minus signs, respectively, $t'_{n,j,\pm}$ denotes the intersection of the horizontal line $y = \varphi_{1,n}(t_{j,+})$ and the line connecting the points $(j\pi/\lambda, 0)$ and $(t_{j,+}, 1)$ [$(j\pi/\lambda, 0)$ and $(t_{j,-}, 1)$]. Therefore,

$$t'_{n,j,\pm} = j\frac{\pi}{\lambda} \pm C_1^* n^{-\zeta/\mu} \left(j + \frac{1}{2}\right)^{(\rho+\nu)/\mu},$$

where $C_1^* > 0$.

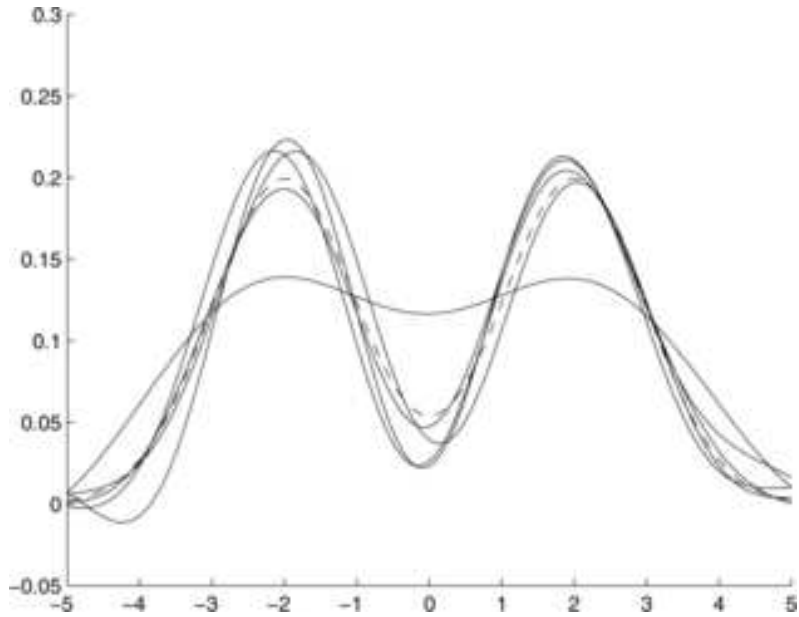

Fig. 7. $n = 400$, $r = 2$, $\rho = 2$, $\xi$ by CV average integrated squared error $= 0.003$.

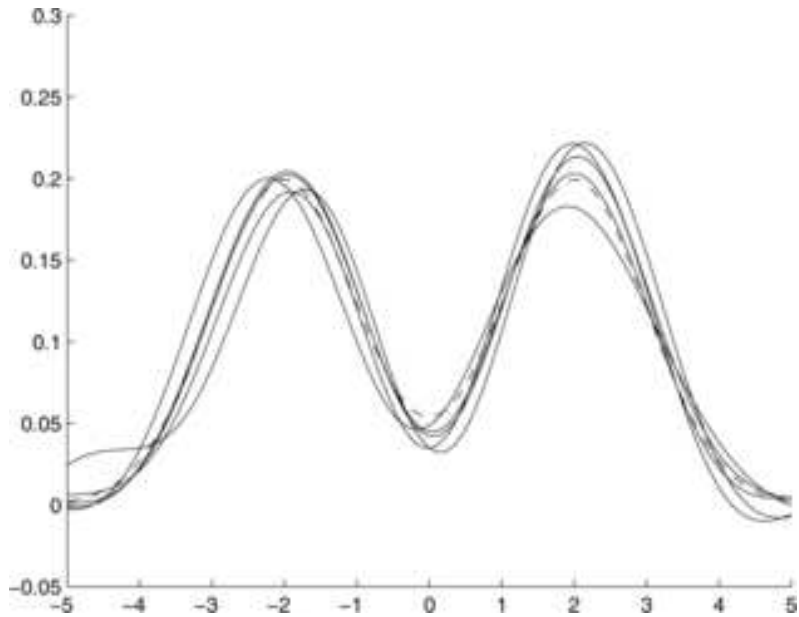

Fig. 8. $n = 700$, $r = 2$, $\rho = 2$, $\xi$ by CV average integrated squared error $= 0.001$.

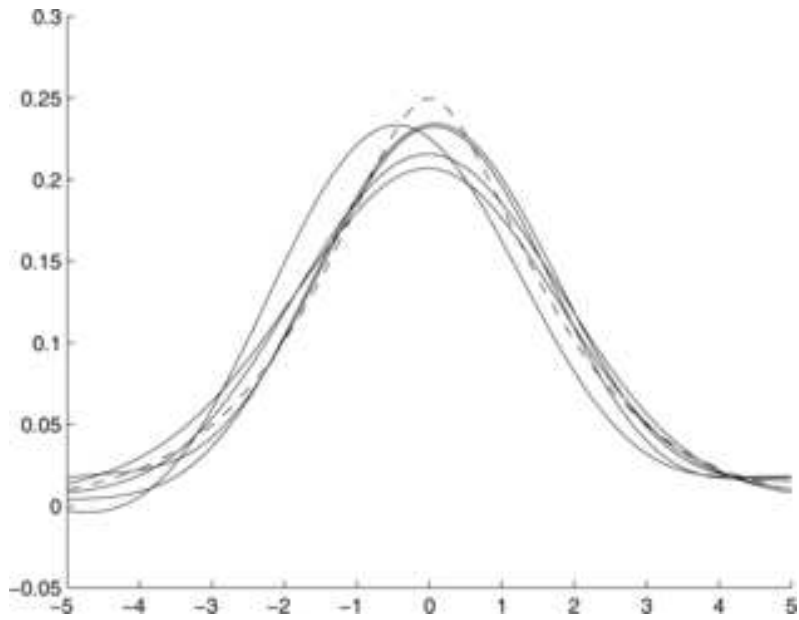

FIG. 9. $n = 400$, $r = 2$, $\rho = 2$, $\xi$ by CV average integrated squared error $= 0.002$.

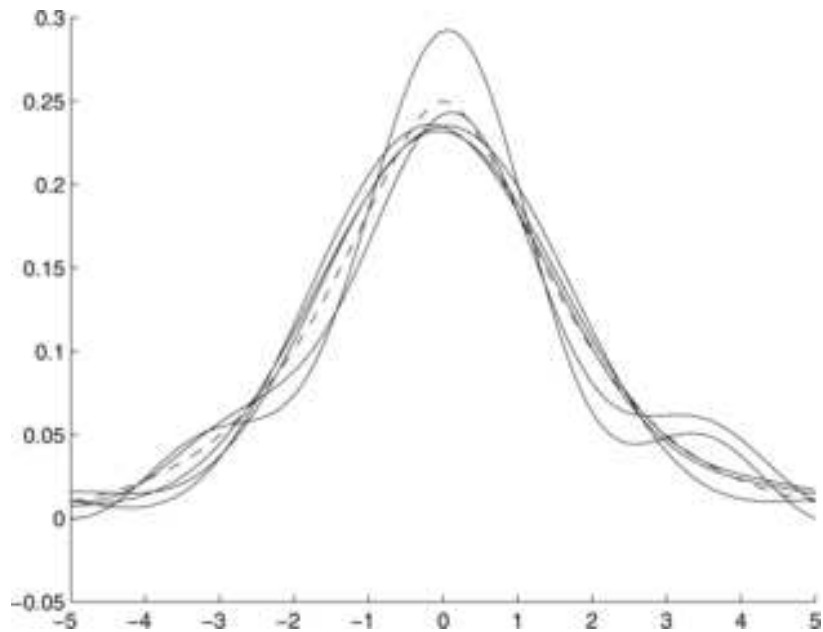

FIG. 10. $n = 700$, $r = 2$, $\rho = 2$, $\xi$ by CV average integrated squared error $= 0.001$.

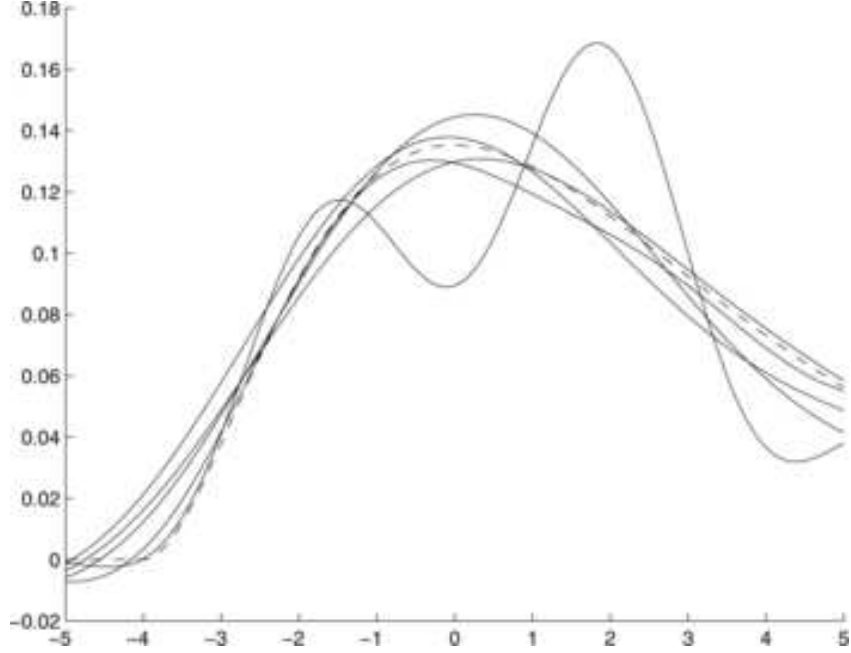

FIG. 11. $n = 400$, $r = 2$, $\rho = 2$, $\xi$ by CV average integrated squared error $= 0.001$.

In the arguments immediately below, leading to a bound on $V_{1,n}$, we use the property $G_{n,j} \subseteq [t'_{n,j,-}, t'_{n,j,+}]$ for $j \leq j'_n \sim n^{\zeta/(\rho+\nu)}$, where $j'_n$ denotes the largest $j$ for which the inclusion relation holds. For larger values of $j$, we use the property $G_{n,j} \subseteq I_j$, thus obtaining

$$V_{1,n} \leq O(n^{2\zeta(r+2)-1}) \left[ \sum_{j=1}^{j'_n} \int_{t'_{n,j,-}}^{t'_{n,j,+}} \{\sin(\lambda t)\}^{2\mu(r+1)} |t|^{-2(\nu+2\rho)-2r(\nu+\rho)} \, dt \right.$$

$$\left. + \int_{t > t'_{j'_n,-}} \{\sin(\lambda t)\}^{2\mu(r+1)} |t|^{-2(\nu+2\rho)-2r(\nu+\rho)} \, dt \right]$$

$$\leq O(n^{2\zeta(r+2)-1}) \sum_{j=1}^{j'_n} j^{-2(\nu+2\rho)-2r(\nu+\rho)} \int_{t'_{n,j,-} - j\pi/\lambda}^{t'_{n,j,+} - j\pi/\lambda} t^{2\mu(r+1)} \, dt$$

$$+ O(n^{2\zeta(r+2)-1} j'^{1-2(\nu+2\rho)-2r(\nu+\rho)}_n)$$

$$\leq O(n^{-1+\zeta(2\mu-1)/\mu}) \sum_{j=1}^{j'_n} j^{-2\rho+(\rho+\nu)/\mu} + O(n^{-1+\zeta(2\nu+1)/(\nu+\rho)}).$$

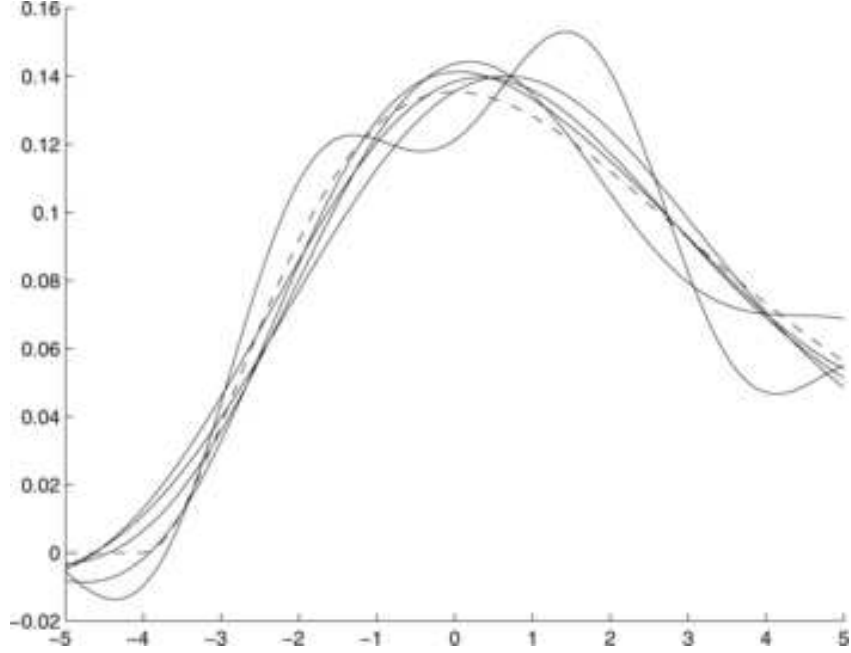

Fig. 12. $n = 700$, $r = 2$, $\rho = 2$, $\xi$ by $CV$ average integrated squared error $= 0.0007$.

From this bound it may be proved that

$$V_{1,n} = \begin{cases} O(n^{-1+\zeta(2\mu-1)/\mu}), & \text{if } \rho > (\mu+\nu)/(2\mu-1), \\ O(n^{-1+\zeta(2\mu-1)/\mu} \log n), & \text{if } \rho = (\mu+\nu)/(2\mu-1), \\ O(n^{-1+\zeta(2\nu+1)/(\nu+\rho)}), & \text{otherwise.} \end{cases}$$

A similar argument can be used to bound the bias term $B_n$:

$$B_n \le \frac{1}{2\pi} \sup_{f_X} \sum_{j=1}^{j'_n} \int_{t'_{n,j,-}}^{t'_{n,j,+}} |f_X^{\text{ft}}(t)|^2 \, dt + O(j'^{-2\beta}_n).$$

At this point it is necessary to treat separately the cases $f_X \in \mathcal{F}^j_{\beta C}$, for $j = 2, 3$. When $j = 3$,

$$B_n \le \frac{1}{2\pi} \sup_{f_X} \sum_{j=1}^{j'_n} j^{-2\beta-1}(t'_{n,j,-} - t'_{n,j,+}) + O(n^{-2\zeta\beta/(\rho+\nu)})$$

$$= \begin{cases} O(n^{-\zeta/\mu}), & \text{if } \rho + \nu < 2\mu\beta, \\ O(n^{-\zeta/\mu} \log n), & \text{if } \rho + \nu = 2\mu\beta, \\ O(n^{-2\beta\zeta/(\rho+\nu)}), & \text{otherwise.} \end{cases}$$

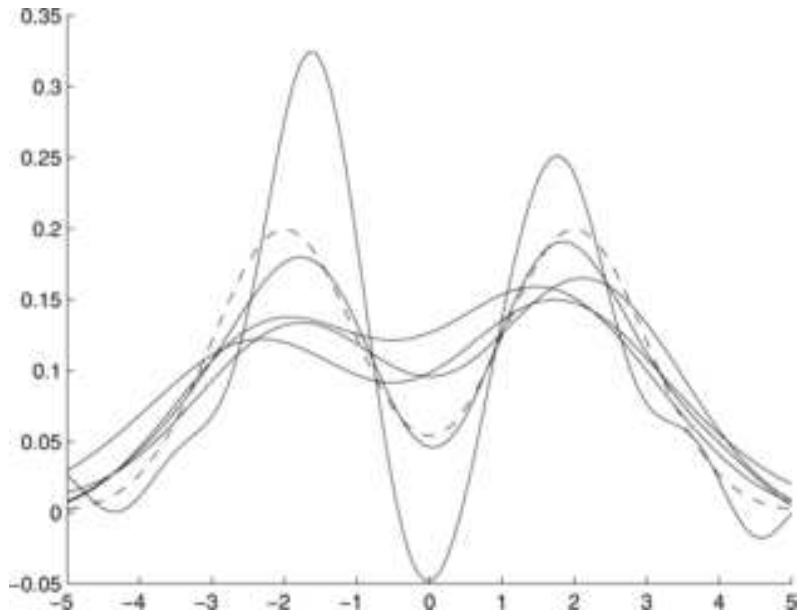

Fig. 13. $n = 400$, bandwidth by CV average integrated squared error $= 0.0012$.

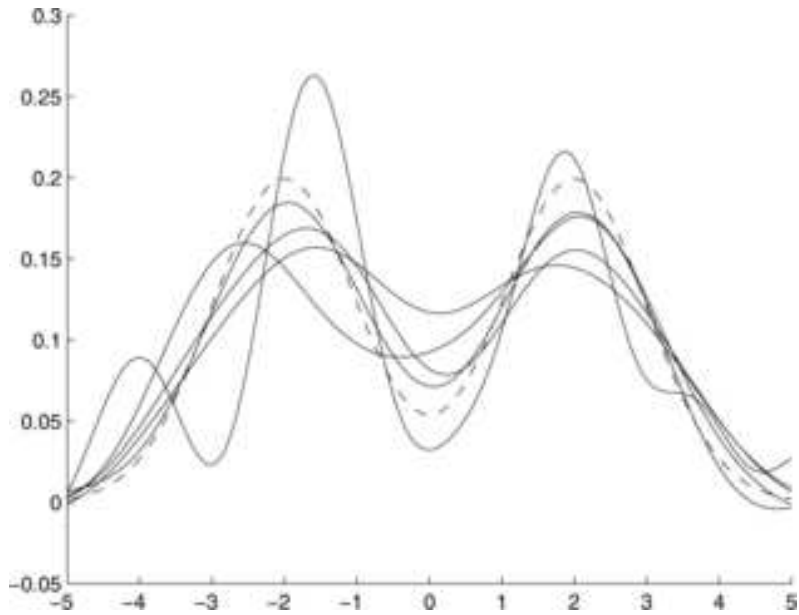

Fig. 14. $n = 700$, bandwidth by CV average integrated squared error $= 0.009$.

When $j = 2$,

$$\max_{t \in [t_{n,j,-}, t_{n,j,+}]} |f_X^{\text{ft}}(t)|^2 = \left\{ \max_{t \in [t_{n,j,-}, t_{n,j,+}]} |f_X^{\text{ft}}(t)|^2 - \min_{t \in [t_{n,j,-}, t_{n,j,+}]} |f_X^{\text{ft}}(t)|^2 \right\}$$
$$+ \min_{t \in [t_{n,j,-}, t_{n,j,+}]} |f_X^{\text{ft}}(t)|^2$$
$$\leq O(j^{-2\beta}) \int_{t'_{n,j,-}}^{t'_{n,j,+}} (1+t^2)^\beta \{|f_X^{\text{ft}}(t)|^2 + |f_X^{\text{ft}}(t)||f_X^{\text{ft}\prime}(t)|\} \, dt.$$

This leads to the same upper bound on $B_n$ as in the case $j = 3$.

Finally we bound $V_{2,n}$. There, we need a lower bound for $G_n$, obtained as follows. Then $t \in G_{n,j}$ is implied by

$$|\sin(\lambda t)| \leq C_2^{-1/\mu} n^{-\zeta/\mu} |t|^{(\rho+\nu)/\mu} = \varphi_{2,n}(t),$$

say. Let $t''_{n,j,-} < t''_{n,j,+}$ denote the intersections of the horizontal line $y = \varphi_{2,n}(t_{j,-})$ and both tangent lines at $t = j\pi/\lambda$ of the curve with equation $y = |\sin(\lambda t)|$. Then, for sufficiently small $j$, $[t''_{n,j,-}, t''_{n,j,+}] \subseteq G_{n,j}$, and so

$$t''_{n,j,\pm} = j \frac{\pi}{\lambda} \pm C_2^* n^{-\zeta/\mu} \left( j - \frac{1}{2} \right)^{(\rho+\nu)/\mu},$$

where $C_2^* > 0$.

Let $j''_n \sim n^{\zeta/(\rho+\nu)}$ denote the largest $j$ for which $[t''_{n,j,-}, t''_{n,j,+}] \subseteq G_{n,j}$. We use this inclusion relation if $j \leq j''_n$, and the relation $I_j \subseteq G_{n,j}$ otherwise. This leads to the property

$$[0, \infty) \cap G_n^c \subseteq \bigcup_{j=1}^{j''_n} (I_j \setminus [t''_{n,j,-}, t''_{n,j,+}]) \cup I_0.$$

To bound $V_{2,n}$ we apply this formula to the integral over $t \in G_n^c$ with $t > 0$. The contribution from $I_0$ can be shown to equal $O(n^{-1})$, and so is negligible. Defining $I''_{n,j} = I_j \setminus [t''_{n,j,-}, t''_{n,j,+}]$, and the shifted set $I^{\text{sh}}_{n,j} = [-\pi/(2\lambda), t''_{n,j,-} - j\pi/\lambda] \cup [t''_{n,j,+} - j\pi/\lambda, \pi/(2\lambda)]$, it may be proved that

$$V_{2,n} \leq O(n^{-1}) \sum_{j=1}^{j''_n} \int_{I''_{n,j}} |\sin(\lambda t)|^{-2\mu} |t|^{2\nu} \, dt$$
$$\leq O(n^{-1}) \sum_{j=1}^{j''_n} j^{2\nu} \int_{I^{\text{sh}}_{n,j}} |\sin(\lambda t)|^{-2\mu} \, dt$$
$$= O(n^{-1+\zeta(2\mu-1)/\mu}) \sum_{j=1}^{j''_n} j^{2\nu+(\rho+\nu)(1-2\mu)/\mu}.$$

This leads to the same upper bound for $V_{2,n}$ that we derived earlier for $V_{1,n}$. Substituting for $\rho$ and $\zeta$ from (4.6) and (4.7), we obtain Theorem 4.2(a).

6.3. *Proof of Theorem* 4.2(b). Using Lemma 6.1, it may be shown that $V_n = O(n^{-1+4\zeta})$. Note too that, for $t \in G_n$,

(6.1) $$|\sin(\lambda t)| \leq C_1^{-1/\mu} n^{-\zeta/\mu} \exp\{(d/\mu)|t|^\gamma\}.$$

Define $t_n = \{\zeta/(2d)\}^{1/\gamma}(\log n)^{1/\gamma}$. If $|t| \leq t_n$, then $|\sin(\lambda t)| \leq C_1^{-1/\mu} n^{-\zeta/(2\mu)}$, and hence, using (6.1), $|t - k\pi/\lambda| \leq \{\pi C^{-1/\mu}/(2\lambda)\} n^{-\zeta/(2\mu)}$ for integers $k$. From these properties it may be proved that, for $f_X \in \mathcal{F}_{\beta C}^j$ for $j = 1, 2$ or $3$,

$$B_n = O(t_n n^{-\zeta/(2\mu)} + t_n^{-2\beta}) = O\{(\log n)^{-2\beta/\gamma}\}.$$

6.4. *Proof of Theorem* 4.3. First we derive (4.10). We introduce the density $f_0(x) = (1 - \cos x)/(\pi x^2)$, having the "tent"-shaped Fourier transform, $f_0^{\text{ft}}(t) = 1 - |t|$ for $|t| \leq 1$, and the supersmooth Cauchy density $f_1(x) = (1/\pi)(1+x^2)^{-1}$. Define too the densities

$$f_{n,\theta}(x) = \tfrac{1}{2}\varepsilon_n f_1(\varepsilon_n x) + \tfrac{1}{2}\varepsilon_n f_0(\varepsilon_n x)\{1 + \eta\theta\cos(Tx)\},$$

with $\eta \in (0, \tfrac{1}{2}]$, $\theta \in \{0, 1\}$ and a positive-valued sequence $\varepsilon_n \downarrow 0$, to be defined later. The characteristic function corresponding to $f_{n,\theta}$ is

$$f_{n,\theta}^{\text{ft}}(x) = \tfrac{1}{2}f_1^{\text{ft}}(t/\varepsilon_n) + \tfrac{1}{2}f_0^{\text{ft}}(t/\varepsilon_n) + \tfrac{1}{4}\eta\theta f_0^{\text{ft}}\{(t+T)/\varepsilon_n\} + \tfrac{1}{4}\eta\theta f_0^{\text{ft}}\{(t-T)/\varepsilon_n\}.$$

Using the fact that $\varepsilon_n \downarrow 0$, it can be shown that if $\eta$ is sufficiently small, then $f_{n,\theta} \in \mathcal{F}_{\beta C}^j$ can be verified for any $C$ and $\beta$.

Write $f_1 * f_2$ for the convolution of functions $f_1$ and $f_2$. It was proved by Fan [12] that if the $\chi^2$-distance between the densities $f_{n,0} * f_\delta$ and $f_{n,1} * f_\delta$ satisfies

(6.2) $$\chi^2(f_{n,1} * f_\delta, f_{n,0} * f_\delta) = \int \frac{|(f_{n,0} * f_\delta)(x) - (f_{n,1} * f_\delta)(x)|^2}{(f_{n,0} * f_\delta)(x)} \, dx = O(n^{-1}),$$

then, for a constant $c > 0$, for any estimator $\hat{f}$ and for all sufficiently large $n$,

(6.3) $$\sup_{f \in \mathcal{F}_{\beta C}^j} E\|\hat{f} - f\|^2 \geq c\|f_{n,1} - f_{n,0}\|^2.$$

It may be shown from the definition of $f_{n,0}(x)$ that

$$\chi^2(f_{n,1} * f_\delta, f_{n,0} * f_\delta) \leq 2\varepsilon_n^{-1} \int \frac{[\{(f_{n,0} - f_{n,1}) * f_\delta\}(x)]^2}{\{f_1(\varepsilon_n \cdot) * f_\delta\}(x)} \, dx.$$

Also, if $q$ is so large so that $\int_{|y| \leq q} f_\delta(y) \, dy > 0$,

$$\{f_1(\varepsilon_n \cdot) * f_\delta\}(x) \geq \pi^{-1} \int_{|y| \leq q} f_\delta(y)\{1 + 2\varepsilon_n^2(x^2 + y^2)\}^{-1} \, dy$$

$$\geq \text{const.}(1 + \varepsilon_n^2 x^2)^{-1}.$$

Therefore,

$$\chi^2(f_{n,1} * f_\delta, f_{n,0} * f_\delta) \leq \text{const.} \int \{(f_{n,0} - f_{n,1}) * f_\delta\}^2(x)(\varepsilon_n^{-1} + \varepsilon_n x^2)\, dx.$$

To bound the right-hand side, use the identity $(f^{\text{ft}})' = i\{\cdot f(\cdot)\}^{\text{ft}}$ and the Parseval identity to obtain

$$\chi^2(f_{n,1} * f_\delta, f_{n,0} * f_\delta) \leq \text{const.} \bigg[\varepsilon_n^{-1} \int |(f_{n,0} - f_{n,1})^{\text{ft}}(t)|^2 |f_\delta^{\text{ft}}(t)|^2\, dt$$

$$+ \varepsilon_n \int |\{(f_{n,0} - f_{n,1})^{\text{ft}}\}'(t)|^2 |f_\delta^{\text{ft}}(t)|^2\, dt$$

$$+ \varepsilon_n \int |(f_{n,0} - f_{n,1})^{\text{ft}}(t)|^2 |(f_\delta^{\text{ft}})'(t)|^2\, dt \bigg].$$

Therefore, using the fact that $(f_{n,0} - f_{n,1})^{\text{ft}}$ and $\{(f_{n,0} - f_{n,1})^{\text{ft}}\}'$ are both supported on $[-T-\varepsilon_n, -T+\varepsilon_n] \cup [T-\varepsilon_n, T+\varepsilon_n]$, we may show that $\chi^2(f_{n,1} * f_\delta, f_{n,0} * f_\delta)$ equals

$$O\bigg\{\varepsilon_n^{-1} \int_{T-\varepsilon_n}^{T+\varepsilon_n} |f_\delta^{\text{ft}}(t)|^2\, dt + \varepsilon_n \int_{T-\varepsilon_n}^{T+\varepsilon_n} |(f_\delta^{\text{ft}})'(t)|^2\, dt + \varepsilon_n^{-1} \int_{T-\varepsilon_n}^{T+\varepsilon_n} |f_\delta^{\text{ft}}(t)|^2\, dt\bigg\}.$$

Using (4.9) to bound the right-hand side, we may thus show that $\chi^2(f_{n,1} * f_\delta, f_{n,0} * f_\delta) = O(\varepsilon_n^{2\mu})$. Hence, choosing $\varepsilon_n = n^{-1/(2\mu)}$ guarantees the validity of (6.3). Again, Parseval's identity may be used to prove that

$$\|f_{n,1} - f_{n,0}\|^2 \geq \text{const.}\varepsilon_n,$$

which implies (4.10).

Finally we turn to (4.11). With the densities $f_0, f_1$ as above, we construct the subclass of densities

$$\tilde{f}_{n,\theta}(x) = \tfrac{1}{2}\{f_1(x) + f_0(x)\} + \text{const.} \sum_{k_n \leq j \leq 2k_n} \theta_j j^{-\beta-(1/2)} \cos(2jx) f_0(x),$$

with $k_n$ an integer satisfying $k_n \uparrow \infty$, and $\theta_j \in \{0, 1\}$. The class $\mathcal{F}_{\beta C}^3$ contains all densities $\tilde{f}_{n,\theta}$. By modeling the $\theta_j$'s as independent random variables with $P(\theta_j = 0) = \tfrac{1}{2}$, the lower bound (4.11) may be established using arguments similar to those in the proof of Fan [15].

DEPARTMENT OF MATHEMATICS AND STATISTICS
THE UNIVERSITY OF MELBOURNE
MELBOURNE, VICTORIA 3010
AUSTRALIA
E-MAIL: halpstat@ms.unimelb.edu.au

INSTITUT FÜR STOCHASTIK UND ANWENDUNGEN
UNIVERSITÄT STUTTGART
PFAFFENWALDRING 57
D-70569 STUTTGART
GERMANY
E-MAIL: meistear@mathematik.uni-stuttgart.de